\documentclass[a4paper,11pt,twoside,reqno]{amsart}

\usepackage[utf8]{inputenc}
\usepackage[plainpages=false,pdfpagelabels=true]{hyperref}
\usepackage{amssymb,amsthm}
\usepackage[margin=1in]{geometry}

\newtheorem{Theorem}{Theorem}[section]
\newtheorem{Prop}[Theorem]{Proposition}
\newtheorem{Lem}[Theorem]{Lemma}

\newtheorem{Cor}[Theorem]{Corollary}

\theoremstyle{definition}

\newtheorem{Bem}[Theorem]{Remark}

\newcommand{\tr}{\operatorname{Tr}}
\newcommand{\dv}{\text{ }dV}
\allowdisplaybreaks[1]

\renewcommand{\epsilon}{\varepsilon}

\newcommand{\R}{\ensuremath{\mathbb{R}}}
\newcommand{\N}{\ensuremath{\mathbb{N}}}

\numberwithin{equation}{section}

\title{A structure theorem for polyharmonic maps between Riemannian manifolds}

\author{Volker Branding}

\address{University of Vienna, Faculty of Mathematics\\
Oskar-Morgenstern-Platz 1, 1090 Vienna, Austria}
\email{volker.branding@univie.ac.at}

\date{\today}



\begin{document}

\begin{abstract}
We prove that polyharmonic maps of arbitrary order from complete nonparabolic Riemannian manifolds
to arbitrary Riemannian manifolds must be harmonic if certain smallness and integrability conditions hold.
\end{abstract}

\maketitle

\section{Introduction and results}
Finding interesting maps between Riemannian manifolds is one of the core problems 
in the field of geometric variational problems. In order to find such maps, one usually studies
a certain energy functional that one associates to the map.
One may then calculate the critical points of this energy functional
yielding the corresponding Euler-Lagrange equation. Very often, 
the solutions of the Euler-Lagrange equation contain interesting geometric data.

For a map \(\phi\colon M\to N\) between two Riemannian manifolds \((M,g)\) and \((N,h)\) much effort has been paid in studying the critical points 
of their energy which is defined by 
\begin{align}
\label{1-energy}
E(\phi)=E_1(\phi)=\int_M|d\phi|^2\dv.
\end{align}
The critical points of \eqref{1-energy} are governed by the vanishing of the so-called 
\emph{tension field} that is
\begin{align}
\label{harmonic-map-equation}
0=\tr_g\bar\nabla d\phi.
\end{align}
Here, we denote the connection on \(\phi^\ast TN\) by \(\bar\nabla\).
Solutions of \eqref{harmonic-map-equation} are called \emph{harmonic maps}. 
The harmonic map equation is a second order semilinear elliptic partial differential equation.

Another energy functional that receives growing attention in the mathematical literature is the 
\emph{bienergy} of a map which is given by
\begin{align*}
E_2(\phi)=\int_M|\tau(\phi)|^2\dv.
\end{align*}
Its critical points are of fourth order and are characterized by the vanishing
of the \emph{bitension field} that is 
\begin{align}
\label{biharmonic-map-equation}
0=\tau_2(\phi):=\bar\Delta\tau(\phi)-\tr_g R^N(d\phi,\tau(\phi))d\phi.
\end{align}
Here, \(\bar\Delta\) denotes the connection Laplacian on \(\phi^\ast TN\) and by \(R^N\) we denote
the curvature tensor on \(N\).
Solutions of \eqref{biharmonic-map-equation} are called \emph{biharmonic maps}.
In contrast to the harmonic map equation, the equation for biharmonic maps is of fourth order
such that its analysis comes with additional technical difficulties.

One natural generalization of both harmonic and biharmonic maps can be obtained
by studying the critical points of the \emph{\(k\)-energy} of a map between two Riemannian manifolds
which is defined by
\begin{align}
\label{k-energy}
E_k(\phi)=\int_M|(d+\delta)^k\phi|^2\dv.
\end{align}
The study of this energy was proposed by Eells and Lemaire in 1983 \cite[p.77, Problem (8.7)]{MR703510}.

The first variation of \eqref{k-energy} was calculated in \cite{MR3007953},
its critical points are referred to as \emph{polyharmonic maps (of order \(k\))} or \(k\)-harmonic maps.
We have to distinguish two cases:
\begin{enumerate}
 \item If \(k=2s,s=1,2,\ldots\), the critical points of \eqref{k-energy} are given by
\begin{align}
\label{tension-2s}
0=\tau_{2s}(\phi):=&\bar\Delta^{2s-1}\tau(\phi)-R^N(\bar\Delta^{2s-2}\tau(\phi),d\phi(e_j))d\phi(e_j) \\
\nonumber&-\sum_{l=1}^{s-1}\bigg(R^N(\bar\nabla_{e_j}\bar\Delta^{s+l-2}\tau(\phi),\bar\Delta^{s-l-1}\tau(\phi))d\phi(e_j) \\
\nonumber&\hspace{1cm}-R^N(\bar\Delta^{s+l-2}\tau(\phi),\bar\nabla_{e_j}\bar\Delta^{s-l-1}\tau(\phi))d\phi(e_j)
\bigg).
\end{align}
\item If \(k=2s+1,s=0,1,\ldots\), the critical points of \eqref{k-energy} are given by
\begin{align}
\label{tension-2s+1}
0=\tau_{2s+1}(\phi):=&\bar\Delta^{2s}\tau(\phi)-R^N(\bar\Delta^{2s-1}\tau(\phi),d\phi(e_j))d\phi(e_j)\\
\nonumber&-\sum_{l=1}^{s-1}\bigg(R^N(\bar\nabla_{e_j}\bar\Delta^{s+l-1}\tau(\phi),\bar\Delta^{s-l-1}\tau(\phi))d\phi(e_j) \\
\nonumber&-R^N(\bar\Delta^{s+l-1}\tau(\phi),\bar\nabla_{e_j}\bar\Delta^{s-l-1}\tau(\phi))d\phi(e_j)
\bigg) \\
&\nonumber-R^N(\bar\nabla_{e_j}\bar\Delta^{s-1}\tau(\phi),\bar\Delta^{s-1}\tau(\phi))d\phi(e_j).
\end{align}
\end{enumerate}
Here, we have set \(\bar\Delta^{-1}=0\), \(\{e_j\},j=1,\ldots,\dim M\) denotes an orthonormal basis of \(TM\)
and we are applying the Einstein summation convention.

The second variation of the \(k\)-energy \eqref{k-energy} has also been calculated in \cite{MR3007953}.

There exists a vast number of results on existence and qualitative behavior of
harmonic maps (\(k=1\)) and biharmonic maps (\(k=2\)).
The polyharmonic map equation is a semilinear elliptic partial differential equation of order \(2k\)
and so far polyharmonic maps did not receive a lot of attention as their governing equation
comes with a large number of derivatives.

Regarding the regularity of weak solutions, some results for polyharmonic maps could be obtained,
see for example \cite{MR2520907, MR2542730}.
However, most of the articles dealing with such kind of questions assume that \(N\) is realized as a submanifold
of Euclidean space. In this case the \(k\)-energy and its critical points depend also on the embedding
in the ambient space. Such kind of maps are usually referred to as \emph{extrinsic polyharmonic maps}.

It follows directly from \eqref{tension-2s} and \eqref{tension-2s+1} that a harmonic map
always solves the equation for polyharmonic maps. This is to be expected as we may rewrite the \(k\)-energy
for a map as follows
\begin{align*}
E_k(\phi)=\int_M|(d+\delta)^{k-2}\delta d\phi|^2\dv=\int_M|(d+\delta)^{k-2}\tau(\phi)|^2\dv
\end{align*}
and thus a harmonic map is always an absolute minimum of the \(k\)-energy for \(k\geq 2\).
However, a \(k\)-harmonic map does not necessarily have to be harmonic, see for example 
\cite{MR3711937} for proper \(k\)-harmonic immersions into spheres.

Let us mention several results on polyharmonic maps that are connected to the main theorem of this article.
In \cite{MR3371364} the authors show that triharmonic immersions 
into a manifold of non-positive constant curvature are minimal under certain
integrability assumptions.

In \cite{MR3809656} the authors show that polyharmonic maps from complete non-compact Riemannian manifolds
to Euclidean space are harmonic under certain integrability assumptions. 
This result has been generalized in \cite{MR3314128}.
Several results on polyharmonic submanifolds of Euclidean space can be found in \cite{MR2911957}.

Polyharmonic maps in space forms have been studied in \cite{MR2869168}
and it is shown that for \(k\geq 3\) a \(k\)-harmonic isometric immersion into a Riemannian manifold 
of nonzero constant sectional curvature is minimal.

The main result of this article shows that polyharmonic maps between Riemannian
manifolds must be harmonic if one assumes certain smallness and integrability conditions
generalizing a previous result on biharmonic maps \cite{MR3834926,branding2018nonexistence}
to polyharmonic maps of arbitrary order. Our result is quite general in the sense that we do not have
to impose any curvature condition on the target and are able to treat polyharmonic maps
of arbitrary high order \(k\).

\begin{Theorem}
\label{main-theorem}
Let \((M,g)\) be a complete non-compact Riemannian manifold 
that admits an Euclidean type Sobolev inequality
and let \(\phi\colon M\to N\) be a polyharmonic map of order \(k\).

For \(k\geq 2\) being even and \(n>2k-2\), the following classification result holds true:
\begin{enumerate}
 \item Suppose that the following condition holds
\begin{align}
\label{assumption-main-a}
\int_M|\bar\nabla^qd\phi|^\frac{n}{q+1}\dv<\epsilon,\qquad \textrm{ for all } ~~0\leq q\leq k-2
\end{align}
 for some \(\epsilon>0\) small enough.
 \item In addition, assume that
 \begin{align}
 \label{assumption-main-b}
 \sum_{q=0}^{k-2}\int_M|\bar\nabla^q\bar\Delta^{\frac{k}{2}-1}\tau(\phi)|^2\dv<\infty.
 \end{align}
\end{enumerate}
Then \(\phi\) must be harmonic.
\par\medskip
For \(k\geq 3\) being odd and \(n>2k\), the following classification result holds true:
\begin{enumerate}
 \item Suppose that the following condition holds
\begin{align}
\label{odd-assumption-main-a}
\int_M|\bar\nabla^qd\phi|^\frac{n}{q+1}\dv<\epsilon,\qquad \textrm{ for all } ~~0\leq q\leq k-1
\end{align}
 for some \(\epsilon>0\) small enough.
 \item In addition, assume that
 \begin{align}
 \label{odd-assumption-main-b}
 \sum_{q=0}^{k-1}\int_M|\bar\nabla^q\bar\Delta^{\frac{k-3}{2}}\tau(\phi)|^2\dv<\infty.
 \end{align}
\end{enumerate}
Then \(\phi\) must be harmonic.
\end{Theorem}

\begin{Bem}
\begin{enumerate}
\item In non-technical terms Theorem \ref{main-theorem} states the following: A polyharmonic map of order \(k\)
 has to be harmonic if its derivatives of order \(1\) up to order \(k-1\) are small and the derivatives of order
 \(k\) up to order \(2k-2\) are bounded. 
\item Note that the assumptions \eqref{assumption-main-a} and \eqref{odd-assumption-main-a} in Theorem \ref{main-theorem} are
natural in the sense that the required smallness conditions
are invariant under rescalings of the metric on the domain.
\end{enumerate}
\end{Bem}

If we also demand positive Ricci curvature on the domain then we obtain
the following variant of Theorem \ref{main-theorem}.

\begin{Cor}
Let \(\phi\colon M\to N\) be a polyharmonic map of order \(k\) and suppose
that the assumptions of Theorem \ref{main-theorem} hold.
If \(M\) has positive Ricci curvature then \(\phi\) must be trivial.
\end{Cor}
This corollary follows from Theorem \ref{main-theorem} together with \cite[Theorem 1.5]{MR1333821}.

In the proof of our main result we will apply an \emph{Euclidean type Sobolev inequality}
of the following form
\begin{align}
\label{sobolev-inequality}
(\int_M|u|^{2n/(n-2)}\dv)^\frac{n-2}{n}\leq C_2\int_M|Du|^2\dv
\end{align}
for all \(u\in W^{1,2}(M)\) with compact support
where \(C_2\) is a positive constant that depends on the geometry of \(M\).
Such an inequality holds in \(\R^n\) and is well-known as \emph{Gagliardo-Nirenberg inequality}.
However, if one considers a complete non-compact Riemannian manifold of infinite volume
one has to make additional assumptions to ensure that an equality of the form \eqref{sobolev-inequality} 
holds. In technical terms a complete Riemannian manifold is called \emph{nonparabolic}
if it admits an inequality of the form \eqref{sobolev-inequality}.

For more details under which assumptions an inequality of the form \eqref{sobolev-inequality}
holds we refer to the introduction of \cite{MR3886921} and references therein.

Throughout this article we will use the following sign conventions:
For the Riemannian curvature
tensor we use \(R(X,Y)Z=[\nabla_X,\nabla_Y]Z-\nabla_{[X,Y]}Z\).
By \(\bar\nabla\) we denote the connection on \(\phi^\ast TN\) and 
for the (rough) Laplacian on \(\phi^\ast TN\) we use
\(\bar\Delta:=-\tr_g(\bar\nabla\bar\nabla-\bar\nabla_\nabla)\).

We will employ Latin letters for indices on the domain \(M\) ranging from \(1\) to \(n\)
and we frequently use the Einstein summation convention, i.e. we will always sum over repeated indices.
We will use the notation \(\{e_j\},j=1,\ldots,n\) to denote an orthonormal basis of \(TM\).

\section{Proof of the main result}
In order to prove the main result Theorem \ref{main-theorem} we have to distinguish between polyharmonic
maps of even and odd order.

Throughout the proof we will make use of a cutoff function  \(0\leq\eta\leq 1\) on \(M\) that satisfies
\begin{align*}
\eta(x)=1\textrm{ for } x\in B_R(x_0),\qquad \eta(x)=0\textrm{ for } x\in B_{2R}(x_0),\qquad |D^q\eta|\leq\frac{C}{R^q}\textrm{ for } x\in M,
\end{align*}
where \(B_R(x_0)\) denotes the geodesic ball around the point \(x_0\) with radius \(R\) and \(q\in\N\).

In addition, we will use the same symbol \(\epsilon\) to denote various quantities that are 
small without referring to their explicit expression.

By iterating \eqref{sobolev-inequality} we obtain
\begin{align}
\label{sobolev-inequality-k}
(\int_M|u|^\frac{2n}{n-2r}\dv)^\frac{n-2r}{2n}\leq C_r(\int_M|D^ru|^2\dv)^\frac{1}{2},\qquad n>2r,
\end{align}
which we will frequently use throughout the proof.
Here, \(C_r\) is a positive constant that depends on \(n,r\) and \(C_2\).
We will use the symbol \(C\) to denote a generic positive constant whose
value may change from line to line.

Note that we will always set \(\bar\nabla^{-1}=\bar\Delta^{-1}=0\).

\subsection{The even case}
In this section we prove Theorem \ref{main-theorem} in the case that we have a polyharmonic
map of even order that is a solution of \eqref{tension-2s}. Hence, we set \(k=2s,s=1,2,\ldots\).
\begin{Lem}
Let \(\phi\colon M\to N\) be a smooth solution of \eqref{tension-2s}
and suppose that \(M\) admits an Euclidean type Sobolev inequality.
In addition, assume that \(n>4s-2\).
Then the following estimate holds
\begin{align}
\label{estimate-a}
\int_M\eta^2|\bar\nabla&\bar\Delta^{2s-2}\tau(\phi)|^2\dv\\
\nonumber\leq 
&C|D\eta|^2\int_M|\bar\Delta^{2s-2}\tau(\phi)|^2\dv
+C_2\big(\int_M|d\phi|^n\dv\big)^\frac{2}{n}\int_M\big|D(\eta|\bar\Delta^{2s-2}\tau(\phi)|)\big|^2\dv \\
\nonumber&
+C\sum_{l=1}^{s-1}\bigg(\big(\int_M|d\phi|^n\dv\big)^\frac{1}{n}
\big(\int_M|\bar\Delta^{s-l-1}\tau(\phi)|^\frac{n}{2(s-l)}\dv\big)^\frac{2(s-l)}{n}\\
\nonumber&\times 
\big(\int_M\big|D^{2(s-l)}(\eta|\bar\nabla\bar\Delta^{s+l-2}\tau(\phi)|)\big|^2\dv\big)^\frac{1}{2}
\big(\int_M\big|D(\eta|\bar\Delta^{2s-2}\tau(\phi)|)\big|^2\dv\big)^\frac{1}{2}\bigg) \\
\nonumber&+ C\sum_{l=1}^{s-1}\bigg(\big(\int_M|d\phi|^n\dv\big)^\frac{1}{n}
\big(\int_M|\bar\nabla\bar\Delta^{s-l-1}\tau(\phi)|^\frac{n}{2(s-l)+1}\dv\big)^\frac{2(s-l)+1}{n}\\
\nonumber&\times 
\big(\int_M\big|D^{2(s-l)+1}(\eta|\bar\Delta^{s+l-2}\tau(\phi)|)\big|^2\dv\big)^\frac{1}{2}
\big(\int_M\big|D(\eta|\bar\Delta^{2s-2}\tau(\phi)|)\big|^2\dv\big)^\frac{1}{2}\bigg),
\end{align}
where the positive constant \(C\) depends on \(n,s\) and the geometries of \(M\) and \(N\).
\end{Lem}
\begin{proof}
We test \eqref{tension-2s} with \(\eta^2\bar\Delta^{2s-2}\tau(\phi)\) and obtain
\begin{align*}
\int_M\eta^2\langle\bar\Delta^{2s-1}\tau(\phi),&\bar\Delta^{2s-2}\tau(\phi)\rangle\dv=\\
&\int_M\eta^2\langle R^N(\bar\Delta^{2s-2}\tau(\phi),d\phi(e_j))d\phi(e_j),\bar\Delta^{2s-2}\tau(\phi)\rangle\dv \\
&+\sum_{l=1}^{s-1}\int_M\eta^2\langle R^N(\bar\nabla_{e_j}\bar\Delta^{s+l-2}\tau(\phi),\bar\Delta^{s-l-1}\tau(\phi))d\phi(e_j),\bar\Delta^{2s-2}\tau(\phi)\rangle\dv \\
&-\sum_{l=1}^{s-1}\int_M\eta^2\langle R^N(\bar\Delta^{s+l-2}\tau(\phi),\bar\nabla_{e_j}\bar\Delta^{s-l-1}\tau(\phi))d\phi(e_j),\bar\Delta^{2s-2}\tau(\phi)\rangle\dv.
\end{align*}
Using integration by parts this yields the following inequality
\begin{align*}
\int_M\eta^2|\bar\nabla\bar\Delta^{2s-2}\tau(\phi)|^2\dv\leq& 2\int_M\eta D\eta\langle\bar\nabla\bar\Delta^{2s-2}\tau(\phi),\bar\Delta^{2s-2}\tau(\phi)\rangle\dv \\
&+C\int_M\eta^2|\bar\Delta^{2s-2}\tau(\phi)|^2|d\phi|^2\dv\\
&+C\sum_{l=1}^{s-1}\int_M|\eta^2\bar\nabla\bar\Delta^{s+l-2}\tau(\phi)||\bar\Delta^{s-l-1}\tau(\phi)||d\phi||\bar\Delta^{2s-2}\tau(\phi)|\dv \\
&+C\sum_{l=1}^{s-1}\int_M|\eta^2\bar\Delta^{s+l-2}\tau(\phi)||\bar\nabla\bar\Delta^{s-l-1}\tau(\phi)||d\phi||\bar\Delta^{2s-2}\tau(\phi)|\dv.
\end{align*}
As a next step we estimate all the terms on the right hand side.
First of all, we find
\begin{align*}
2\int_M\eta D\eta\langle\bar\nabla\bar\Delta^{2s-2}\tau(\phi),\bar\Delta^{2s-2}\tau(\phi)\rangle\dv
\leq& \frac{1}{2}\int_M\eta^2|\bar\nabla\bar\Delta^{2s-2}\tau(\phi)|^2\dv \\
&+C|D\eta|^2\int_M|\bar\Delta^{2s-2}\tau(\phi)|^2\dv,
\end{align*}
where we applied Young's inequality.

As a next step we calculate
\begin{align*}
\int_M\eta^2|\bar\Delta^{2s-2}\tau(\phi)|^2|d\phi|^2\dv&\leq \big(\int_M|d\phi|^n\dv\big)^\frac{2}{n}\big(\int_M(\eta|\bar\Delta^{2s-2}\tau(\phi)|)^\frac{2n}{n-2}\dv\big)^\frac{n-2}{n} \\
&\leq C_2\big(\int_M|d\phi|^n\dv\big)^\frac{2}{n}\int_M\big|D(\eta|\bar\Delta^{2s-2}\tau(\phi)|)\big|^2\dv,
\end{align*}
where we first used Hölder's inequality and applied \eqref{sobolev-inequality} afterwards
assuming that \(n>2\).

Moreover, we find
\begin{align*}
\int_M&\eta^2|\bar\nabla\bar\Delta^{s+l-2}\tau(\phi)||\bar\Delta^{s-l-1}\tau(\phi)||d\phi||\bar\Delta^{2s-2}\tau(\phi)|\dv \\
\leq
&\big(\int_M|d\phi|^n\dv\big)^\frac{1}{n}
\big(\int_M|\bar\Delta^{s-l-1}\tau(\phi)|^\frac{n}{2(s-l)}\dv\big)^\frac{2(s-l)}{n}\\
&\times\big(\int_M(\eta|\bar\nabla\bar\Delta^{s+l-2}\tau(\phi)|)^\frac{2n}{n-4s+4l}\dv\big)^\frac{n-4s+4l}{2n}
\big(\int_M(\eta|\bar\Delta^{2s-2}\tau(\phi)|)^\frac{2n}{n-2}\dv\big)^\frac{n-2}{2n}\\
\leq & C\big(\int_M|d\phi|^n\dv\big)^\frac{1}{n}
\big(\int_M|\bar\Delta^{s-l-1}\tau(\phi)|^\frac{n}{2(s-l)}\dv\big)^\frac{2(s-l)}{n}\\
&\times 
\big(\int_M\big|D^{2(s-l)}(\eta|\bar\nabla\bar\Delta^{s+l-2}\tau(\phi)|)\big|^2\dv\big)^\frac{1}{2}
\big(\int_M\big|D(\eta|\bar\Delta^{2s-2}\tau(\phi)|)\big|^2\dv\big)^\frac{1}{2},
\end{align*}
where we applied \eqref{sobolev-inequality-k} and thus have to impose that \(n>4s-4\).

Finally, we compute
\begin{align*}
\int_M&\eta^2|\bar\Delta^{s+l-2}\tau(\phi)||\bar\nabla\bar\Delta^{s-l-1}\tau(\phi)||d\phi||\bar\Delta^{2s-2}\tau(\phi)|\dv \\
\leq
&\big(\int_M|d\phi|^n\dv\big)^\frac{1}{n}
\big(\int_M|\bar\nabla\bar\Delta^{s-l-1}\tau(\phi)|^\frac{n}{2(s-l)+1}\dv\big)^\frac{2(s-l)+1}{n}\\
&\times\big(\int_M(\eta|\bar\Delta^{s+l-2}\tau(\phi)|)^\frac{2n}{n-4s+4l-2}\dv\big)^\frac{n-4s+4l-2}{2n}
\big(\int_M(\eta|\bar\Delta^{2s-2}\tau(\phi)|)^\frac{2n}{n-2}\dv\big)^\frac{n-2}{2n}\\
\leq & C\big(\int_M|d\phi|^n\dv\big)^\frac{1}{n}
\big(\int_M|\bar\nabla\bar\Delta^{s-l-1}\tau(\phi)|^\frac{n}{2(s-l)+1}\dv\big)^\frac{2(s-l)+1}{n}\\
&\times 
\big(\int_M\big|D^{2(s-l)+1}(\eta|\bar\Delta^{s+l-2}\tau(\phi)|)\big|^2\dv\big)^\frac{1}{2}
\big(\int_M\big|D(\eta|\bar\Delta^{2s-2}\tau(\phi)|)\big|^2\dv\big)^\frac{1}{2},
\end{align*}
where we again applied \eqref{sobolev-inequality-k} under the condition that \(n>4s-2\).
The claim follows by combining the estimates.
\end{proof}

\begin{Lem}
Let \(\phi\colon M\to N\) be a smooth solution of \eqref{tension-2s}
and suppose that \(M\) admits an Euclidean type Sobolev inequality.
In addition, assume that \(n>4s-2\).
If the \(n\)-energy of \(\phi\) is sufficiently small, that is
\begin{align}
\label{assumption-b}
\int_M|d\phi|^n\dv<\epsilon           
\end{align} for some small \(\epsilon>0\),
then the following inequality holds
\begin{align}
\label{estimate-b}
&(1-\epsilon C_2)\int_M\eta^2|\bar\nabla\bar\Delta^{2s-2}\tau(\phi)|^2\dv \\
\nonumber\leq & \frac{C}{R^2}\int_M|\bar\Delta^{2s-2}\tau(\phi)|^2\dv \\
\nonumber&+C\sum_{l=1}^{s-1}\epsilon
\big(
\frac{C}{R^2}\int_M|\bar\Delta^{2s-2}\tau(\phi)|^2\dv
+\int_M\eta^2|\bar\nabla\bar\Delta^{2s-2}\tau(\phi)|^2\dv
\big)^\frac{1}{2} \\
\nonumber&
\hspace{0.4cm}\times\bigg(\big(\int_M|\bar\Delta^{s-l-1}\tau(\phi)|^\frac{n}{2(s-l)}\dv\big)^\frac{2(s-l)}{n}\\
\nonumber&\hspace{0.6cm}\times 
\big(
\int_M\eta^2|\bar\nabla^{2(s-l)+1}\bar\Delta^{s+l-2}\tau(\phi)|^2\dv
+\sum_{q=1}^{2(s-l)}\frac{1}{R^{2q}}\int_M|\bar\nabla^{2(s-l)-q}\bar\nabla\bar\Delta^{s+l-2}\tau(\phi)|^2\dv
\big)^\frac{1}{2} \\
\nonumber&\hspace{0.6cm}+ 
\big(
\int_M|\bar\nabla\bar\Delta^{s-l-1}\tau(\phi)|^\frac{n}{2(s-l)+1}\dv\big)^\frac{2(s-l)+1}{n}\\
\nonumber&\hspace{0.6cm}\times
\big(
\int_M\eta^2|\bar\nabla^{2(s-l)+1}\bar\Delta^{s+l-2}\tau(\phi)|^2 \dv
+\sum_{q=1}^{2(s-l)+1}\frac{1}{R^{2q}}\int_M|\bar\nabla^{2(s-l)-q}\bar\nabla\bar\Delta^{s+l-2}\tau(\phi)|^2\dv
\big)^\frac{1}{2}\bigg).
\end{align}
\end{Lem}

\begin{proof}
First of all, we note that
\begin{align*}
\int_M\big|D(\eta|\bar\Delta^{2s-2}\tau(\phi)|)\big|^2\dv\leq
\frac{C}{R^2}\int_M|\bar\Delta^{2s-2}\tau(\phi)|^2\dv
+\int_M\eta^2|\bar\nabla\bar\Delta^{2s-2}\tau(\phi)|^2\dv,
\end{align*}
where we used the properties of the cutoff function \(\eta\).
In addition, we find
\begin{align*}
\big|D^{2(s-l)}(\eta|\bar\nabla\bar\Delta^{s+l-2}\tau(\phi)|)\big|^2=&\big|\sum_{q=0}^{2(s-l)}C_qD^q\eta|\bar\nabla^{2(s-l)-q}\bar\nabla\bar\Delta^{s+l-2}\tau(\phi)|\big|^2\\
\leq& C\sum_{q=0}^{2(s-l)}|D^q\eta|^2||\bar\nabla^{2(s-l)-q}\bar\nabla\bar\Delta^{s+l-2}\tau(\phi)|^2\\
=&C\eta^2|\bar\nabla^{2(s-l)+1}\bar\Delta^{s+l-2}\tau(\phi)|^2 \\
&+C\sum_{q=1}^{2(s-l)}\frac{1}{R^{2q}}|\bar\nabla^{2(s-l)-q}\bar\nabla\bar\Delta^{s+l-2}\tau(\phi)|^2
\end{align*}
and similarly for the other term of the same structure.

Applying the estimate \eqref{estimate-a} and making use of \eqref{assumption-b} completes the proof.
\end{proof}

\begin{Lem}
\label{lem-c}
Let \(\phi\colon M\to N\) be a smooth solution of \eqref{tension-2s}
and suppose that \(M\) admits an Euclidean type Sobolev inequality.
In addition, assume that \(n>4s-2\).
Moreover, suppose that the \(n\)-energy of \(\phi\) is sufficiently small, that is
\begin{align*}
\int_M|d\phi|^n\dv<\epsilon           
\end{align*} for some small \(\epsilon>0\).
In addition, we assume that
\begin{align} 
 \label{assumption-c2}
 \int_M|\bar\Delta^{r}\tau(\phi)|^\frac{n}{2(r+1)}\dv&<\infty\qquad \text{ for all } ~~ 0\leq r \leq s-2, \\
 \nonumber \int_M|\bar\nabla\bar\Delta^{r}\tau(\phi)|^\frac{n}{2(r+1)+1}\dv&<\infty\qquad \text{ for all } ~~ 0\leq r \leq s-2 
 \end{align}
and 
\begin{align}
\label{assumption-c3}
\int_M|\bar\nabla^{r}\bar\Delta^{s-2}\tau(\phi)|^2\dv<\infty\qquad \text{ for all } ~~ 2\leq r \leq 2s.
\end{align}
Then the following inequality holds
\begin{align}
\label{estimate-c}
\int_M\eta^2|&\bar\nabla\bar\Delta^{2s-2}\tau(\phi)|^2\dv 
\leq C\epsilon\sum_{l=1}^{s-1}
\int_M\eta^2|\bar\nabla^{2(s-l)+1}\bar\Delta^{s+l-2}\tau(\phi)|^2\dv,
\end{align}
where the positive constant \(C\) depends on \(n,s,\epsilon\) and the geometries of \(M\) and \(N\).
\end{Lem}

\begin{proof}
Taking the limit \(R\to\infty\) in \eqref{estimate-b} and using the integrability assumption \eqref{assumption-c3} yields
\begin{align*}
(1-&\epsilon C_2)\int_M\eta^2|\bar\nabla\bar\Delta^{2s-2}\tau(\phi)|^2\dv \\
\leq 
&C\epsilon\sum_{l=1}^{s-1}
\big(\int_M\eta^2|\bar\nabla\bar\Delta^{2s-2}\tau(\phi)|^2\dv
\big)^\frac{1}{2} \\
&
\hspace{0.5cm}\times\bigg(\big(\int_M|\bar\Delta^{s-l-1}\tau(\phi)|^\frac{n}{2(s-l)}\dv\big)^\frac{2(s-l)}{n}
\big(
\int_M\eta^2|\bar\nabla^{2(s-l)+1}\bar\Delta^{s+l-2}\tau(\phi)|^2\dv
\big)^\frac{1}{2} \\
&\hspace{0.7cm}+ 
\big(
\int_M|\bar\nabla\bar\Delta^{s-l-1}\tau(\phi)|^\frac{n}{2(s-l)+1}\dv\big)^\frac{2(s-l)+1}{n}
\big(
\int_M\eta^2|\bar\nabla^{2(s-l)+1}\bar\Delta^{s+l-2}\tau(\phi)|^2 \dv
\big)^\frac{1}{2}\bigg).
\end{align*}
After employing the assumptions \eqref{assumption-c2} we are left with
\begin{align*}
\int_M\eta^2|&\bar\nabla\bar\Delta^{2s-2}\tau(\phi)|^2\dv \\
&\leq C\epsilon
\big(\int_M\eta^2|\bar\nabla\bar\Delta^{2s-2}\tau(\phi)|^2\dv
\big)^\frac{1}{2} 
\sum_{l=1}^{s-1}
\big(
\int_M\eta^2|\bar\nabla^{2(s-l)+1}\bar\Delta^{s+l-2}\tau(\phi)|^2\dv
\big)^\frac{1}{2},
\end{align*}
which completes the proof.
\end{proof}

\begin{Lem}
Suppose that the assumptions of Lemma \ref{lem-c} hold true.
Then we have the following inequality
\begin{align}
\label{estimate-d}
\int_M\eta^2|&\bar\nabla\bar\Delta^{2s-2}\tau(\phi)|^2\dv 
\leq C\epsilon
\int_M\eta^2|\bar\nabla^{2s-1}\bar\Delta^{s-1}\tau(\phi)|^2\dv.
\end{align}
\end{Lem}
\begin{proof}
As we can always estimate \(|\bar\Delta\tau(\phi)|\leq\sqrt{n}|\bar\nabla^2\tau(\phi)|\) 
the statement follows from \eqref{estimate-c}.
\end{proof}

The estimate \eqref{estimate-d} is already almost the estimate that we are looking for.
Our strategy to conclude that \(\bar\nabla\bar\Delta^{2s-2}\tau(\phi)=0\) will now be the following: 
We will interchange as many derivatives on the right-hand side of \eqref{estimate-d} as needed in order 
to obtain a suitable power of the Laplacian. Making use of the smallness assumption on \(\epsilon\)
we can then absorb the right-hand side into the left. However, before we can do so, we have to control
all curvature terms that appear when turning the covariant derivatives into powers of the Laplacian.

\begin{Lem}
Let \(\phi\colon M\to N\) be a smooth map.
Then the following inequality holds
\begin{align}
\label{estimate-e}
\frac{1}{2}\int_M\eta^2|&\bar\nabla^{2s-1}\bar\Delta^{s-1}\tau(\phi)|^2\dv \\
\nonumber\leq & C\int_M\eta^2|\bar\nabla\bar\Delta^{2s-2}\tau(\phi)|^2\dv 
+\frac{C}{R^2}\sum_{r=0}^{s-2}\int_M|\bar\nabla^{2s-2-2r}\bar\Delta^{s-1+r}\tau(\phi)|^2\dv\\
&\nonumber+\sum_{r=0}^{s-2}\int_M\eta^2|\bar\nabla^{2s-2-2r}\bar\Delta^{s-1+r}\tau(\phi)||[\bar\Delta,\bar\nabla^{2s-2-2r}]\bar\Delta^{s-1+r}\tau(\phi)|\dv.
\end{align}
\end{Lem}
\begin{proof}
Using integration by parts and interchanging covariant derivatives we find
\begin{align*}
\int_M\eta^2|\bar\nabla^{2s-1}\bar\Delta^{s-1}\tau(\phi)|^2\dv=
&-2\int_M\eta D\eta\langle\bar\nabla^{2s-2}\bar\Delta^{s-1}\tau(\phi),\bar\nabla^{2s-1}\bar\Delta^{s-1}\tau(\phi)\rangle\dv \\
&+2\int_M\eta D\eta\langle\bar\nabla^{2s-2}\bar\Delta^{s-1}\tau(\phi),\bar\nabla^{2s-3}\bar\Delta^{s}\tau(\phi)\rangle\dv \\
&-\int_M\eta^2\langle\bar\nabla^{2s-2}\bar\Delta^{s-1}\tau(\phi),[\bar\Delta,\bar\nabla^{2s-2}]\bar\Delta^{s-1}\tau(\phi)\rangle\dv \\
&+\int_M\eta^2\langle\bar\nabla^{2s-3}\bar\Delta^s\tau(\phi),\bar\Delta\bar\nabla^{2s-3}\bar\Delta^{s-1}\tau(\phi)\rangle\dv.
\end{align*}
Using that 
\begin{align*}
|\langle\bar\nabla^{2s-3}\bar\Delta^s\tau(\phi),\bar\Delta\bar\nabla^{2s-3}\bar\Delta^{s-1}\tau(\phi)\rangle|\leq 
C|\bar\nabla^{2s-3}\bar\Delta^s\tau(\phi)|^2+\frac{1}{2}|\bar\nabla^{2s-1}\bar\Delta^{s-1}\tau(\phi)|^2
\end{align*}
and also applying Young's inequality to the terms involving the derivative of the cutoff function \(\eta\) and estimating 
the commutator term we find
\begin{align*}
\frac{1}{2}\int_M\eta^2|\bar\nabla^{2s-1}\bar\Delta^{s-1}\tau(\phi)|^2\dv\leq&C\int_M\eta^2|\bar\nabla^{2s-3}\bar\Delta^{s}\tau(\phi)|^2\dv 
+\frac{C}{R^2}\int_M|\bar\nabla^{2s-2}\bar\Delta^{s-1}\tau(\phi)|^2\dv\\
&+\int_M\eta^2|\bar\nabla^{2s-2}\bar\Delta^{s-1}\tau(\phi)||[\bar\Delta,\bar\nabla^{2s-2}]\bar\Delta^{s-1}\tau(\phi)|\dv.
\end{align*}
The claim now follows by iterating the above procedure \((s-1)\) times.
\end{proof}

In the following we will often employ the so-called \(\star\)-notation.
Here, a \(\star\) refers to various contractions between the objects involved.
The following Lemma is a variant of the calculations performed in \cite[Section 3]{MR3830277}.

\begin{Lem}
Let \(\phi\colon M\to N\) be a smooth map.
Then the following identity holds
\begin{align}
\label{commutator}
[\bar\Delta,\bar\nabla^r]X=\sum_{\sum_{l_i}+\sum_{m_j}=r}\nabla^{l_1}R^N\star \underbrace{\bar\nabla^{m_1}d\phi\star\ldots\bar\nabla^{m_{l_1}}d\phi}_{l_1-\text{times}}
\star\bar\nabla^{l_2}d\phi\star\bar\nabla^{l_3}d\phi\star\bar\nabla^{l_4}X,
\end{align}
where \(r\in\N\) and \(X\in\Gamma(\phi^\ast TN)\).
\end{Lem}
\begin{proof}
The formula holds for \(r=1\) by the following explicit calculation
\begin{align*}
[\bar\Delta,\bar\nabla_i]X=&(\nabla_{d\phi(e_j)}R^N)(d\phi(e_j),d\phi(e_i))X
+R^N(\tau(\phi),d\phi(e_i))X \\
&+R^N(d\phi(e_j),\bar\nabla_j(d\phi(e_i)))X
+2R^N(d\phi(e_j),d\phi(e_i))\bar\nabla_jX.
\end{align*}
The claim follows by iteration.
\end{proof}

\begin{Lem}
Let \(\phi\colon M\to N\) be a smooth map and
suppose that \(M\) admits an Euclidean type Sobolev inequality.
In addition, assume that \(n>4s-2\). Then the following estimate holds
\begin{align}
\label{estimate-f1}
\sum_{r=0}^{s-2}&\int_M\eta^2|\bar\nabla^{2s-2-2r}\bar\Delta^{s-1+r}\tau(\phi)||[\bar\Delta,\bar\nabla^{2s-2-2r}]\bar\Delta^{s-1+r}\tau(\phi)|\dv \\
\nonumber\leq &
C\big(\frac{1}{R^2}\int_M|\bar\nabla^{2s-2}\bar\Delta^{s-1}\tau(\phi)|^2\dv
+\int_M\eta^2|\bar\nabla^{2s-1}\bar\Delta^{s-1}\tau(\phi)|^2\dv\big)^\frac{1}{2}\\
\nonumber&\times\sum_{r=0}^{s-1}\sum_{\sum_{l_i}+\sum_{m_j}=2s-2-2r}
\big(\int_M|\bar\nabla^{l_2}d\phi|^\frac{n}{l_2+1}\dv\big)^\frac{l_2+1}{n}
\big(\int_M|\bar\nabla^{l_3}d\phi|^\frac{n}{l_3+1}\dv\big)^\frac{l_3+1}{n}\\
\nonumber&\times\big(
\int_M\eta^2|\bar\nabla^{2s-1}\bar\Delta^{s-1}\tau(\phi)|^2\dv
+\sum_{q=1}^{2s-1-l_4-2r}\frac{C}{R^{2q}}\int_M|\bar\nabla^{2s-1-q}\bar\Delta^{s-1}\tau(\phi)|^2\dv\big)^\frac{1}{2}\\
&\nonumber\times\prod_{j=1}^{l_1}\big(\int_M|\bar\nabla^{m_j}d\phi|^\frac{n}{m_j+l_1}\dv\big)^\frac{m_j+l_1}{n}.
\end{align}
\end{Lem}

\begin{proof}
Using \eqref{commutator} we find 
\begin{align*}
\int_M&\eta^2|\bar\nabla^{2s-2-2r}\bar\Delta^{s-1+r}\tau(\phi)||[\bar\Delta,\bar\nabla^{2s-2-2r}]\bar\Delta^{s-1+r}\tau(\phi)|\dv \\
&=\sum_{\sum_{l_i}+\sum_{m_j}=2s-2-2r}\int_M\eta^2|\bar\nabla^{2s-2-2r}\bar\Delta^{s-1+r}\tau(\phi)|
\nabla^{l_1}R^N\star \underbrace{\bar\nabla^{m_1}d\phi\star\ldots\bar\nabla^{m_{l_1}}d\phi}_{l_1-\text{times}}\\
&\hspace{4cm}\star\bar\nabla^{l_2}d\phi\star\bar\nabla^{l_3}d\phi\star\bar\nabla^{l_4}\bar\Delta^{s-1+r}\tau(\phi)|\dv\\
&\leq C\sum_{\sum_{l_i}+\sum_{m_j}=2s-2-2r}\int_M\eta^2|\bar\nabla^{2s-2}\bar\Delta^{s-1}\tau(\phi)|
|\underbrace{\bar\nabla^{m_1}d\phi\star\ldots\bar\nabla^{m_{l_1}}d\phi}_{l_1-\text{times}}|\\
&\hspace{4cm}\times|\bar\nabla^{l_2}d\phi||\bar\nabla^{l_3}d\phi||\bar\nabla^{l_4+2r}\bar\Delta^{s-1}\tau(\phi)|\dv.
\end{align*}
As a next step we calculate
\begin{align*}
&\int_M\eta^2|\bar\nabla^{2s-2}\bar\Delta^{s-1}\tau(\phi)|
|\underbrace{\bar\nabla^{m_1}d\phi\star\ldots\bar\nabla^{m_{l_1}}d\phi}_{l_1-\text{times}}|
|\bar\nabla^{l_2}d\phi||\bar\nabla^{l_3}d\phi||\bar\nabla^{l_4+2r}\bar\Delta^{s-1}\tau(\phi)|\dv \\
\leq&\big(\int_M(\eta|\bar\nabla^{2s-2}\bar\Delta^{s-1}\tau(\phi)|)^\frac{2n}{n-2}\dv\big)^\frac{n-2}{2n}
\big(\int_M|\bar\nabla^{l_2}d\phi|^\frac{n}{l_2+1}\dv\big)^\frac{l_2+1}{n}
\big(\int_M|\bar\nabla^{l_3}d\phi|^\frac{n}{l_3+1}\dv\big)^\frac{l_3+1}{n}\\
&\times\big(\int_M(\eta|\bar\nabla^{l_4+2r}\bar\Delta^{s-1}\tau(\phi)|)^\frac{2n}{n-2(2s-1-l_4-2r)}\dv\big)^\frac{n-2(2s-1-l_4-2r)}{2n}\\
&\times\big(\int_M|\underbrace{\bar\nabla^{m_1}d\phi\star\ldots\bar\nabla^{m_{l_1}}d\phi}_{l_1-\text{times}}|^\frac{n}{\sum_{m_j}+l_1}\dv\big)^\frac{\sum_{m_j}+l_1}{n},
\end{align*}
where we applied Hölder's inequality under the assumption that \(\sum_{l_i}+\sum_{m_j}=2s-2-2r\).

Note that
\begin{align*}
\big(\int_M(\eta|\bar\nabla^{2s-2}\bar\Delta^{s-1}\tau(\phi)|)^\frac{2n}{n-2}\dv\big)^\frac{n-2}{2n}
\leq &C_2\big(\int_M|D(\eta|\bar\nabla^{2s-2}\bar\Delta^{s-1}\tau(\phi)|)|^2\dv\big)^\frac{1}{2} \\
\leq &
\big(\frac{C}{R^2}\int_M|\bar\nabla^{2s-2}\bar\Delta^{s-1}\tau(\phi)|^2\dv \\
&+\int_M\eta^2|\bar\nabla^{2s-1}\bar\Delta^{s-1}\tau(\phi)|^2\dv\big)^\frac{1}{2},
\end{align*}
where we applied \eqref{sobolev-inequality} under the assumption that \(n>2\).

Employing \eqref{sobolev-inequality-k} under the assumption that \(n>4s-2\) we find
\begin{align*}
\big(\int_M&(\eta|\bar\nabla^{l_4+2r}\bar\Delta^{s-1}\tau(\phi)|)^\frac{2n}{n-2(2s-1-l_4-2r)}\dv\big)^\frac{n-2(2s-1-l_4-2r)}{2n} \\
&\leq\big(\int_M\eta^2|\bar\nabla^{2s-1}\bar\Delta^{s-1}\tau(\phi)|^2\dv
+\sum_{q=1}^{2s-1-l_4-2r}\frac{C}{R^{2q}}\int_M|\bar\nabla^{2s-1-q}\bar\Delta^{s-1}\tau(\phi)|^2\dv\big)^\frac{1}{2}.
\end{align*}

In addition, thanks to Hölder's inequality once more we find
\begin{align*}
\big(\int_M|\underbrace{\bar\nabla^{m_1}d\phi\star\ldots\bar\nabla^{m_{l_1}}d\phi}_{l_1-\text{times}}|^\frac{n}{\sum_{m_j}+l_1}\dv\big)^\frac{\sum_{m_j}+l_1}{n}
\leq\prod_{j=1}^{l_1}\big(\int_M|\bar\nabla^{m_j}d\phi|^\frac{n}{m_j+l_1}\dv\big)^\frac{m_j+l_1}{n},
\end{align*}
which yields the claim.
\end{proof}

\begin{Prop}
Let \(\phi\colon M\to N\) be a smooth solution of \eqref{tension-2s} and
suppose that \(M\) admits an Euclidean type Sobolev inequality.
Moreover, suppose that \(n>4s-2\) and assume 
\begin{align}
\label{assumption-identity-a-1}
\int_M|\bar\nabla^qd\phi|^\frac{n}{q+1}\dv<\epsilon,\qquad \textrm{ for all } 0\leq q\leq 2s-2
\end{align}
for some small \(\epsilon>0\).
In addition, assume that
\begin{align}
\label{assumption-identity-a-2}
\sum_{q=0}^{2s-2}\int_M|\bar\nabla^q\bar\Delta^{s-1}\tau(\phi)|^2\dv<\infty.
\end{align}
Then we have
\begin{align}
\label{identity-a}
\bar\nabla\bar\Delta^{2s-2}\tau(\phi)=0.
\end{align}
\end{Prop}
\begin{proof}
Combining \eqref{estimate-e} with the estimates \eqref{estimate-f1} we obtain 
\begin{align*}
\frac{1}{2}\int_M\eta^2|&\bar\nabla^{2s-1}\bar\Delta^{s-1}\tau(\phi)|^2\dv \\
\leq & C\int_M\eta^2|\bar\nabla\bar\Delta^{2s-2}\tau(\phi)|^2\dv 
+\frac{C}{R^2}\sum_{r=0}^{s-2}\int_M|\bar\nabla^{2s-2-2r}\bar\Delta^{s-1+r}\tau(\phi)|^2\dv\\
&+C\big(\frac{1}{R^2}\int_M|\bar\nabla^{2s-2}\bar\Delta^{s-1}\tau(\phi)|^2\dv
+\int_M\eta^2|\bar\nabla^{2s-1}\bar\Delta^{s-1}\tau(\phi)|^2\dv\big)^\frac{1}{2}\\
&\hspace{0.4cm}\times
\sum_{r=0}^{s-2}\sum_{\sum_{l_i}+\sum_{m_j}=2s-2-2r}
\big(\int_M|\bar\nabla^{l_2}d\phi|^\frac{n}{l_2+1}\dv\big)^\frac{l_2+1}{n}
\big(\int_M|\bar\nabla^{l_3}d\phi|^\frac{n}{l_3+1}\dv\big)^\frac{l_3+1}{n}\\
&\hspace{0.4cm}\times\big(
\int_M\eta^2|\bar\nabla^{2s-1}\bar\Delta^{s-1}\tau(\phi)|^2\dv
+\sum_{q=1}^{2s-1-l_4-2r}\frac{C}{R^{2q}}\int_M|\bar\nabla^{2s-1-q}\bar\Delta^{s-1}\tau(\phi)|^2\dv\big)^\frac{1}{2}\\
&\hspace{0.4cm}\times\prod_{j=1}^{l_1}\big(\int_M|\nabla^{m_j}d\phi|^\frac{n}{m_j+l_1}\dv\big)^\frac{m_j+l_1}{n} 
\\
\leq & C\int_M\eta^2|\bar\nabla\bar\Delta^{2s-2}\tau(\phi)|^2\dv 
+\frac{C}{R^2}\int_M|\bar\nabla^{2s-2}\bar\Delta^{s-1}\tau(\phi)|^2\dv\\
&+C\big(\frac{1}{R^2}\int_M|\bar\nabla^{2s-2}\bar\Delta^{s-1}\tau(\phi)|^2\dv
+\int_M\eta^2|\bar\nabla^{2s-1}\bar\Delta^{s-1}\tau(\phi)|^2\dv\big)^\frac{1}{2}\\
&\times\epsilon
\sum_{r=0}^{s-2}\sum_{\sum_{l_i}+\sum_{m_j}=2s-2-2r} \\
&\hspace{0.4cm}\times\big(\int_M\eta^2|\bar\nabla^{2s-1}\bar\Delta^{s-1}\tau(\phi)|^2\dv 
+\sum_{q=1}^{2s-1-l_4-2r}\frac{C}{R^{2q}}\int_M|\bar\nabla^{2s-1-q}\bar\Delta^{s-1}\tau(\phi)|^2\dv\big)^\frac{1}{2},
\end{align*}
where we applied the smallness condition \eqref{assumption-identity-a-1} in the second step.
Taking the limit \(R\to\infty\) while making use of the finiteness assumption \eqref{assumption-identity-a-2}
and using the smallness of \(\epsilon\) we can deduce
\begin{align*}
\int_M\eta^2|\bar\nabla^{2s-1}\bar\Delta^{s-1}\tau(\phi)|^2\dv\leq C\int_M\eta^2|\bar\nabla\bar\Delta^{2s-2}\tau(\phi)|^2\dv.
\end{align*}
Combining this estimate with \eqref{estimate-d} and making use of the smallness of \(\epsilon\) once
more we obtain the claim.
\end{proof}

\begin{Bem}
Note that it was enough to demand the smallness of the \(n\)-energy of \(\phi\)
in order to achieve the estimate \eqref{estimate-d}. However, due to the necessity
of interchanging covariant derivatives we now also had to demand that suitable powers
of higher order derivatives are small \eqref{assumption-identity-a-1}, whereas 
it was enough to demand that they are bounded for the previous estimate \eqref{estimate-c}.
\end{Bem}

In order to prove our main result Theorem \ref{main-theorem} we recall the 
the following result of Gaffney \cite{MR0062490}:
\begin{Theorem}
\label{gaffney}
Let \((M,g)\) be a complete Riemannian manifold. If a \(C^1\) one-form \(\omega\)
satisfies 
\begin{align*}
\int_M|\omega|\dv<\infty \qquad\text{ and }\qquad \int_M|\delta\omega|\dv<\infty
\end{align*}
or, equivalently, a \(C^1\) vector field \(X\) defined by \(\omega(Y)=g(X,Y)\),
satisfies
\begin{align*}
\int_M|X|\dv<\infty \qquad\text{ and }\qquad \int_M\operatorname{div}(X)\dv<\infty,
\end{align*}
then
\begin{align*}
\int_M(\delta\omega)\dv=\int_M\operatorname{div}(X)\dv=0.
\end{align*}
\end{Theorem}

At this point we are ready to complete the proof of Theorem \ref{main-theorem}.
\begin{proof}[Proof of Theorem \ref{main-theorem} (even case)]
We use a suitable iteration method to prove the main result.
\begin{enumerate}
 \item \underline{Step 1}: We define a family of one-forms with \(j\geq 1\) via
\begin{align*}
\omega_j(X)&:=\langle\bar\Delta^j\tau(\phi),\bar\nabla_X\bar\Delta^{j-1}\tau(\phi)\rangle,\qquad s-1\leq j \leq 2s-2.
\end{align*}
If we assume that \(\bar\Delta^j\tau(\phi)\) is parallel, then
\begin{align*}
\delta\omega_j=|\bar\Delta^j\tau(\phi)|^2,\qquad s-1\leq j \leq 2s-2.
\end{align*}
Moreover, we have the estimate
\begin{align*}
\int_M|\omega_j|\dv&\leq\big(\int_M|\bar\Delta^j\tau(\phi)|^2\dv\big)^\frac{1}{2}\big(\int_M|\bar\nabla\bar\Delta^{j-1}\tau(\phi)|^2\dv\big)^\frac{1}{2}.
\end{align*}
Due to Theorem \ref{gaffney} we then find:\\
If \(\bar\Delta^j\tau(\phi)\) is parallel, \(\|\bar\Delta^j\tau(\phi)\|_{L^2}<\infty\) and \(\|\bar\nabla\bar\Delta^{j-1}\tau(\phi)\|_{L^2}<\infty\),
then \(\bar\Delta^j\tau(\phi)=0\).

\item \underline{Step 2}: 
In the second step we assume that \(\bar\Delta^j\tau(\phi)=0\) for \(j\geq 1\) and calculate
\begin{align*}
0=&\int_M\eta^2\langle\bar\Delta^j\tau(\phi),\bar\Delta^{j-1}\tau(\phi)\rangle\dv\\
=&-\int_M\eta^2|\bar\nabla\bar\Delta^{j-1}\tau(\phi)|^2\dv 
-2\int_M\eta D\eta\langle\bar\nabla\bar\Delta^{j-1}\tau(\phi),\bar\Delta^{j-1}\tau(\phi)\rangle\dv.
\end{align*}
We can deduce that
\begin{align*}
\frac{1}{2}\int_M\eta^2|\bar\nabla\bar\Delta^{j-1}\tau(\phi)|^2\dv\leq\frac{C}{R^2}\int_M|\bar\Delta^{j-1}\tau(\phi)|^2\dv ~~\rightarrow 0 \text{ as } R\rightarrow\infty.
\end{align*}
Hence, we may conclude that \(\bar\Delta^{j-1}\tau(\phi)\) is a parallel vector field.

Starting with \(j=2s-2\) using \eqref{identity-a} and after iterating the two steps \((s-1)\) times we find that \(\bar\nabla\bar\Delta^{s-1}\tau(\phi)=0\),
where we made use of the assumptions \eqref{assumption-main-a}.

\item \underline{Step 3}: We define another family of one-forms that fits to the smallness conditions \eqref{assumption-main-b} as follows
\begin{align*}
\omega_r(X):=|\bar\Delta^r\tau(\phi)|^\frac{n-4r-3}{2r+2}\langle\bar\Delta^r\tau(\phi),\bar\nabla_X\bar\Delta^{r-1}\tau(\phi)\rangle,\qquad 1\leq r \leq s-1.
\end{align*}
Note that
\begin{align*}
|\omega_r|\leq C|\bar\nabla^{2r+1}d\phi|^\frac{n-2r-1}{2r+2}|\bar\nabla^{2r}d\phi|
\end{align*}
such that
\begin{align*}
\int_M|\omega_r|\dv\leq C\big(\int_M|\bar\nabla^{2r+1}d\phi|^\frac{n}{2r+2}\dv\big)^\frac{n-2r-1}{n}\big(\int_M|\bar\nabla^{2r}d\phi|^\frac{n}{2r+1}\dv\big)^\frac{2r+1}{n}\leq C
\end{align*}
due to \eqref{assumption-main-b}.

In addition, if we assume that \(\bar\Delta^r\tau(\phi)\) is parallel then we obtain
\begin{align*}
\delta\omega_r=|\bar\Delta^r\tau(\phi)|^\frac{n+1}{2r+2}.
\end{align*}
As \(|\bar\Delta^r\tau(\phi)|\) is constant we find 
\begin{align*}
\int_M|\delta\omega_r|\dv\leq C|\bar\Delta^r\tau(\phi)|^\frac{1}{2r+2}\int_M|\bar\nabla^{2r+1}d\phi|^\frac{n}{2r+2}\dv\leq C
\end{align*}
again due to \eqref{assumption-main-b}.
By Theorem \ref{gaffney} we can then conclude that \(\bar\Delta^r\tau(\phi)=0\).

\item \underline{Step 4}: 
Knowing that \(\bar\Delta^r\tau(\phi)=0\) from the last step we calculate
\begin{align*}
0=&\int_M\eta^2\langle\bar\Delta^r\tau(\phi),\bar\Delta^{r-1}\tau(\phi)\rangle|\bar\Delta^{r-1}\tau(\phi)|^\frac{n-4r}{2r}\dv \\
=&-2\int_M\eta D\eta\langle\bar\nabla\bar\Delta^{r-1}\tau(\phi),\bar\Delta^{r-1}\tau(\phi)\rangle|\bar\Delta^{r-1}\tau(\phi)|^\frac{n-4r}{2r}\dv \\
&-\int_M\eta^2|\bar\nabla\bar\Delta^{r-1}\tau(\phi)|^2|\bar\Delta^{r-1}\tau(\phi)|^\frac{n-4r}{2r}\dv \\
&-\frac{n-4r}{r}\int_M\eta^2|\langle\bar\nabla\bar\Delta^{r-1}\tau(\phi),\bar\Delta^{r-1}\tau(\phi)\rangle|^\frac{n-4r}{2r}\dv.
\end{align*}
Note that \(\frac{n-4r}{r}>0\) due to the assumptions such that
\begin{align*}
\int_M\eta^2|\bar\nabla\bar\Delta^{r-1}\tau(\phi)|^2|\bar\Delta^{r-1}\tau(\phi)|^\frac{n-4r}{2r}\dv
&\leq
\frac{C}{R^2}\int_M|\bar\Delta^{r-1}\tau(\phi)|^\frac{n}{2r}\dv \\
&\leq\frac{C}{R^2}\int_M|\bar\nabla^{2r-1}d\phi|^\frac{n}{2r}\dv ~~\rightarrow 0 \text{ as } R\rightarrow\infty
\end{align*}
again due to \eqref{assumption-main-b}.
Hence, we can conclude that \(\bar\nabla\bar\Delta^{r-1}\tau(\phi)=0\).

After iterating the last two steps for \((s-2)\) times we find that \(\bar\nabla\tau(\phi)=0\).
\item \underline{The final step}: 
Again, we define a one-form \(\omega_0(X)\) via
\begin{align*}
\omega_0(X)&:=|\tau(\phi)|^\frac{n-3}{2}\langle\tau(\phi),d\phi(X)\rangle.
\end{align*}
As \(\tau(\phi)\) is parallel we find
\begin{align*}
\delta\omega_0=|\tau(\phi)|^\frac{n+1}{2}
\end{align*}
and due to the assumptions \eqref{assumption-main-b} we also have
\begin{align*}
\int_M|\delta\omega_0|\dv&\leq C|\tau(\phi)|^\frac{1}{2}\int_M|\nabla d\phi|^\frac{n}{2}\dv\leq C, \\
\int_M|\omega_0|\dv&\leq C\big(\int_M|d\phi|^n\dv\big)^\frac{1}{n}\big(\int_M|\bar\nabla d\phi|^\frac{n}{2}\dv\big)^\frac{n-1}{n}<\infty.
\end{align*}
Applying Theorem \ref{gaffney} we can now conclude that \(\tau(\phi)=0\).
\end{enumerate}
\end{proof}

\subsection{The odd case}
In this section we consider the case of \(k\) being odd that is \(\phi\) is a solution of \eqref{tension-2s+1}.
As the proof is similar to the case of \(k\) being even we do not give as many details as before.
We set \(k=2s+1,s=1,\ldots\) and assume that \(s>0\) as the case \(s=0\) would correspond to \(\phi\) being harmonic.
Note that the assumptions in the odd case of Theorem \ref{main-theorem} are slightly different compared to the even
case.

\begin{Lem}
Let \(\phi\colon M\to N\) be a smooth solution of \eqref{tension-2s+1} and
suppose that \(M\) admits an Euclidean type Sobolev inequality.
Assume that \(n>4s\).
Then the following estimate holds
\begin{align}
\label{odd-estimate-a}
\int_M\eta^2|\bar\nabla&\bar\Delta^{2s-1}\tau(\phi)|^2\dv\\
\nonumber\leq 
&C|D\eta|^2\int_M|\bar\Delta^{2s-1}\tau(\phi)|^2\dv
+C_2\big(\int_M|d\phi|^n\dv\big)^\frac{2}{n}\int_M\big|D(\eta|\bar\Delta^{2s-1}\tau(\phi)|)\big|^2\dv \\
\nonumber&
+C\sum_{l=0}^{s-1}\bigg(\big(\int_M|d\phi|^n\dv\big)^\frac{1}{n}
\big(\int_M|\bar\Delta^{s-l-1}\tau(\phi)|^\frac{n}{2(s-l)}\dv\big)^\frac{2(s-l)}{n}\\
\nonumber&\times 
\big(\int_M\big|D^{2(s-l)}(\eta|\bar\nabla\bar\Delta^{s+l-1}\tau(\phi)|)\big|^2\dv\big)^\frac{1}{2}
\big(\int_M\big|D(\eta|\bar\Delta^{2s-1}\tau(\phi)|)\big|^2\dv\big)^\frac{1}{2}\bigg) \\
\nonumber&+ C\sum_{l=1}^{s-1}\bigg(\big(\int_M|d\phi|^n\dv\big)^\frac{1}{n}
\big(\int_M|\bar\nabla\bar\Delta^{s-l-1}\tau(\phi)|^\frac{n}{2(s-l)+1}\dv\big)^\frac{2(s-l)+1}{n}\\
\nonumber&\times 
\big(\int_M\big|D^{2(s-l)+1}(\eta|\bar\Delta^{s+l-1}\tau(\phi)|)\big|^2\dv\big)^\frac{1}{2}
\big(\int_M\big|D(\eta|\bar\Delta^{2s-1}\tau(\phi)|)\big|^2\dv\big)^\frac{1}{2}\bigg),
\end{align}
where the positive constant \(C\) depends on \(n,s\) and the geometries of \(M\) and \(N\).
\end{Lem}
\begin{proof}
We test \eqref{tension-2s+1} with \(\eta^2\bar\Delta^{2s-1}\tau(\phi)\) and 
after integration by parts this yields the following inequality
\begin{align*}
\frac{1}{2}\int_M\eta^2|\bar\nabla\bar\Delta^{2s-1}\tau(\phi)|^2\dv\leq&C|D\eta|^2\int_M|\bar\Delta^{2s-1}\tau(\phi)|^2\dv\\
&+C\int_M\eta^2|\bar\Delta^{2s-1}\tau(\phi)|^2|d\phi|^2\dv\\
&+C\sum_{l=0}^{s-1}\int_M\eta^2|\bar\nabla\bar\Delta^{s+l-1}\tau(\phi)||\bar\Delta^{s-l-1}\tau(\phi)||d\phi||\bar\Delta^{2s-1}\tau(\phi)|\dv \\
&+C\sum_{l=1}^{s-1}\int_M\eta^2|\bar\Delta^{s+l-1}\tau(\phi)||\bar\nabla\bar\Delta^{s-l-1}\tau(\phi)||d\phi||\bar\Delta^{2s-1}\tau(\phi)|\dv,
\end{align*}
where we used Young's inequality.

As a next step we estimate all the terms on the right hand side starting with
\begin{align*}
\eta^2\int_M|\bar\Delta^{2s-1}\tau(\phi)|^2|d\phi|^2\dv&\leq \big(\int_M|d\phi|^n\dv\big)^\frac{2}{n}\big(\int_M(\eta|\bar\Delta^{2s-1}\tau(\phi)|)^\frac{2n}{n-2}\dv\big)^\frac{n-2}{n} \\
&\leq C_2\big(\int_M|d\phi|^n\dv\big)^\frac{2}{n}\int_M\big|D(\eta|\bar\Delta^{2s-1}\tau(\phi)|)\big|^2\dv,
\end{align*}
where we first used Hölder's inequality and applied \eqref{sobolev-inequality} afterwards.

Moreover, we find
\begin{align*}
\int_M&\eta^2|\bar\nabla\bar\Delta^{s+l-1}\tau(\phi)||\bar\Delta^{s-l-1}\tau(\phi)||d\phi||\bar\Delta^{2s-1}\tau(\phi)|\dv \\
\leq
&\big(\int_M|d\phi|^n\dv\big)^\frac{1}{n}
\big(\int_M|\bar\Delta^{s-l-1}\tau(\phi)|^\frac{n}{2(s-l)}\dv\big)^\frac{2(s-l)}{n}\\
&\times\big(\int_M(\eta|\bar\nabla\bar\Delta^{s+l-1}\tau(\phi)|)^\frac{2n}{n-4s+4l}\dv\big)^\frac{n-4s+4l}{2n}
\big(\int_M(\eta|\bar\Delta^{2s-1}\tau(\phi)|)^\frac{2n}{n-2}\dv\big)^\frac{n-2}{2n}\\
\leq & C\big(\int_M|d\phi|^n\dv\big)^\frac{1}{n}
\big(\int_M|\bar\Delta^{s-l-1}\tau(\phi)|^\frac{n}{2(s-l)}\dv\big)^\frac{2(s-l)}{n}\\
&\times 
\big(\int_M\big|D^{2(s-l)}(\eta|\bar\nabla\bar\Delta^{s+l-1}\tau(\phi)|)\big|^2\dv\big)^\frac{1}{2}
\big(\int_M\big|D(\eta|\bar\Delta^{2s-1}\tau(\phi)|)\big|^2\dv\big)^\frac{1}{2},
\end{align*}
where we applied \eqref{sobolev-inequality-k} and thus have to impose that \(n>4s\).

Finally, we compute
\begin{align*}
\int_M&\eta^2|\bar\Delta^{s+l-1}\tau(\phi)||\bar\nabla\bar\Delta^{s-l-1}\tau(\phi)||d\phi||\bar\Delta^{2s-1}\tau(\phi)|\dv \\
\leq
&\big(\int_M|d\phi|^n\dv\big)^\frac{1}{n}
\big(\int_M|\bar\nabla\bar\Delta^{s-l-1}\tau(\phi)|^\frac{n}{2(s-l)+1}\dv\big)^\frac{2(s-l)+1}{n}\\
&\times\big(\int_M(\eta|\bar\Delta^{s+l-1}\tau(\phi)|)^\frac{2n}{n-4s+4l-2}\dv\big)^\frac{n-4s+4l-2}{2n}
\big(\int_M(\eta|\bar\Delta^{2s-1}\tau(\phi)|)^\frac{2n}{n-2}\dv\big)^\frac{n-2}{2n}\\
\leq & C\big(\int_M|d\phi|^n\dv\big)^\frac{1}{n}
\big(\int_M|\bar\nabla\bar\Delta^{s-l-1}\tau(\phi)|^\frac{n}{2(s-l)+1}\dv\big)^\frac{2(s-l)+1}{n}\\
&\times 
\big(\int_M\big|D^{2(s-l)+1}(\eta|\bar\Delta^{s+l-1}\tau(\phi)|)\big|^2\dv\big)^\frac{1}{2}
\big(\int_M\big|D(\eta|\bar\Delta^{2s-1}\tau(\phi)|)\big|^2\dv\big)^\frac{1}{2},
\end{align*}
where we again applied \eqref{sobolev-inequality-k} under the condition that \(n>4s-2\).
The claim follows by combining the estimates.
\end{proof}

\begin{Lem}
Let \(\phi\colon M\to N\) be a smooth solution of \eqref{tension-2s+1} and
suppose that \(M\) admits an Euclidean type Sobolev inequality.
Assume that \(n>4s\).
Moreover, suppose that the \(n\)-energy of \(\phi\) is sufficiently small, that is
\begin{align}
\label{odd-assumption-b}
\int_M|d\phi|^n\dv<\epsilon           
\end{align} for some small \(\epsilon>0\).
Then the following inequality holds
\begin{align}
\label{odd-estimate-b}
(1-&\epsilon C_2)\int_M\eta^2|\bar\nabla\bar\Delta^{2s-1}\tau(\phi)|^2\dv \\
\nonumber\leq & \frac{C}{R^2}\int_M|\bar\Delta^{2s-1}\tau(\phi)|^2\dv \\
\nonumber &+C\epsilon\sum_{l=0}^{s-1}\bigg(
\big(\int_M|\bar\Delta^{s-l-1}\tau(\phi)|^\frac{n}{2(s-l)}\dv\big)^\frac{2(s-l)}{n}\\
\nonumber&\hspace{0.4cm}\times 
\big(\int_M\eta^2|\bar\nabla^{2(s-l)+1}\bar\Delta^{s+l-1}\tau(\phi)|^2\dv
+\sum_{q=1}^{2(s-l)}\frac{1}{R^{2q}}\int_M|\bar\nabla^{2(s-l)-q}\bar\nabla\bar\Delta^{s+l-1}\tau(\phi)|^2\dv\big)^\frac{1}{2}\\
\nonumber&\hspace{0.4cm}\times
\big(\frac{C}{R^2}\int_M|\bar\Delta^{2s-1}\tau(\phi)|^2\dv
+\int_M\eta^2|\bar\nabla\bar\Delta^{2s-1}\tau(\phi)|^2\dv\big)^\frac{1}{2}\bigg) \\
\nonumber &+C\epsilon\sum_{l=1}^{s-1}\bigg(
\big(\int_M|\bar\nabla\bar\Delta^{s-l-1}\tau(\phi)|^\frac{n}{2(s-l)+1}\dv\big)^\frac{2(s-l)+1}{n}\\
\nonumber&\hspace{0.4cm}\times 
\big(\int_M\eta^2|\bar\nabla^{2(s-l)+1}\bar\Delta^{s+l-1}\tau(\phi)|^2\dv
+\sum_{q=1}^{2(s-l)+1}\frac{1}{R^{2q}}\int_M|\bar\nabla^{2(s-l)+1-q}\bar\Delta^{s+l-1}\tau(\phi)|^2\dv\big)^\frac{1}{2}\\
\nonumber&\hspace{0.4cm}\times
\big(\frac{C}{R^2}\int_M|\bar\Delta^{2s-1}\tau(\phi)|^2\dv
+\int_M\eta^2|\bar\nabla\bar\Delta^{2s-1}\tau(\phi)|^2\dv\big)^\frac{1}{2}\bigg).
\end{align}
\end{Lem}

\begin{proof}
First of all, we note that
\begin{align*}
\int_M\big|D(\eta|\bar\Delta^{2s-1}\tau(\phi)|)\big|^2\dv\leq
\frac{C}{R^2}\int_M|\bar\Delta^{2s-1}\tau(\phi)|^2\dv
+\int_M\eta^2|\bar\nabla\bar\Delta^{2s-1}\tau(\phi)|^2\dv,
\end{align*}
where we used the properties of the cutoff function \(\eta\).
In addition, we find
\begin{align*}
\big|D^{2(s-l)}(\eta|\bar\nabla\bar\Delta^{s+l-1}\tau(\phi)|)\big|^2=&\big|\sum_{q=0}^{2(s-l)}C_qD^q\eta|\bar\nabla^{2(s-l)-q}\bar\nabla\bar\Delta^{s+l-1}\tau(\phi)|\big|^2\\
\leq& C\sum_{q=0}^{2(s-l)}|D^q\eta|^2||\bar\nabla^{2(s-l)-q}\bar\nabla\bar\Delta^{s+l-1}\tau(\phi)|^2\\
=&C\eta^2|\bar\nabla^{2(s-l)+1}\bar\Delta^{s+l-1}\tau(\phi)|^2 \\
&+C\sum_{q=1}^{2(s-l)}\frac{1}{R^{2q}}|\bar\nabla^{2(s-l)-q}\bar\nabla\bar\Delta^{s+l-1}\tau(\phi)|^2
\end{align*}
and similarly for the other term of the same structure.

Applying the estimate \eqref{odd-estimate-a} and making use of \eqref{odd-assumption-b} completes the proof.
\end{proof}

\begin{Lem}
\label{odd-lem-c}
Let \(\phi\colon M\to N\) be a smooth solution of \eqref{tension-2s+1} and
suppose that \(M\) admits an Euclidean type Sobolev inequality.
Assume that \(n>4s\)
and that the \(n\)-energy of \(\phi\) is sufficiently small, that is
\begin{align*}
\int_M|d\phi|^n\dv<\epsilon           
\end{align*}
for some small \(\epsilon>0\).
In addition, we assume that
\begin{align} 
 \label{odd-assumption-d2}
 \int_M|\bar\Delta^{r}\tau(\phi)|^\frac{n}{2r}\dv&<\infty,\qquad \text{ for all } ~~ 0\leq r \leq s-1, \\
 \nonumber \int_M|\bar\nabla\bar\Delta^{r}\tau(\phi)|^\frac{n}{2r+1}\dv&<\infty,\qquad \text{ for all } ~~ 0\leq r \leq s-2 
 \end{align}
and 
\begin{align}
\label{odd-assumption-d3}
\int_M|\bar\nabla^{r}\bar\Delta^{s-1}\tau(\phi)|^2\dv<\infty,\qquad \text{ for all } ~~ 0\leq r \leq 2s-1.
\end{align}
Then the following inequality holds 
\begin{align}
\label{odd-estimate-d}
\int_M\eta^2|\bar\nabla\bar\Delta^{2s-1}\tau(\phi)|^2\dv 
\leq C\epsilon
\int_M\eta^2|\bar\nabla^{2s+1}\bar\Delta^{s-1}\tau(\phi)|^2\dv.
\end{align}
\end{Lem}

\begin{proof}
Taking the limit \(R\to\infty\) in \eqref{odd-estimate-b} and using the integrability assumption \eqref{odd-assumption-d3} yields
\begin{align*}
(1-&\epsilon C_2)\int_M\eta^2|\bar\nabla\bar\Delta^{2s-1}\tau(\phi)|^2\dv \\
\leq 
&C\epsilon
\big(\int_M\eta^2|\bar\nabla\bar\Delta^{2s-1}\tau(\phi)|^2\dv
\big)^\frac{1}{2} \\
&
\hspace{0.5cm}\times
\bigg(\sum_{l=0}^{s-1}\big(\int_M|\bar\Delta^{s-l-1}\tau(\phi)|^\frac{n}{2(s-l)}\dv\big)^\frac{2(s-l)}{n}
\big(
\int_M\eta^2|\bar\nabla^{2(s-l)+1}\bar\Delta^{s+l-1}\tau(\phi)|^2\dv
\big)^\frac{1}{2} \\
&\hspace{0.7cm}+ 
\big(\sum_{l=1}^{s-1}
\int_M|\bar\nabla\bar\Delta^{s-l-1}\tau(\phi)|^\frac{n}{2(s-l)+1}\dv\big)^\frac{2(s-l)+1}{n}
\big(
\int_M\eta^2|\bar\nabla^{2(s-l)+1}\bar\Delta^{s+l-1}\tau(\phi)|^2 \dv
\big)^\frac{1}{2}\bigg).
\end{align*}
After employing the assumptions \eqref{odd-assumption-d2} we are left with
\begin{align*}
\int_M|&\bar\nabla\bar\Delta^{2s-1}\tau(\phi)|^2\eta^2\dv \\
&\leq C\epsilon
\big(\int_M|\bar\nabla\bar\Delta^{2s-1}\tau(\phi)|^2\eta^2\dv
\big)^\frac{1}{2} 
\sum_{l=0}^{s-1}
\big(
\int_M\eta^2|\bar\nabla^{2(s-l)+1}\bar\Delta^{s+l-1}\tau(\phi)|^2\dv
\big)^\frac{1}{2}\\
&\leq C\epsilon
\big(\int_M|\bar\nabla\bar\Delta^{2s-1}\tau(\phi)|^2\eta^2\dv
\big)^\frac{1}{2} 
\big(
\int_M\eta^2|\bar\nabla^{2s+1}\bar\Delta^{s-1}\tau(\phi)|^2\dv
\big)^\frac{1}{2},
\end{align*}
which completes the proof.
\end{proof}

At this point we have to interchange covariant derivatives on the right hand side of \eqref{odd-estimate-d}
similar to the even case. Fortunately, we can use the identities \eqref{commutator}, \eqref{estimate-f1} developed for the even case here as well.

\begin{Prop}
Let \(\phi\colon M\to N\) be a smooth solution of \eqref{tension-2s+1} and
suppose that \(M\) admits an Euclidean type Sobolev inequality.
Moreover, suppose that \(n>4s+2\) and assume 
\begin{align}
\label{odd-assumption-identity-a-1}
\int_M|\bar\nabla^qd\phi|^\frac{n}{q+1}\dv<\epsilon,\qquad \textrm{ for all } 0\leq q\leq 2s
\end{align}
for some small \(\epsilon>0\).
In addition, assume that
\begin{align}
\label{odd-assumption-identity-a-2}
\sum_{q=0}^{2s}\int_M|\bar\nabla^q\bar\Delta^{s-1}\tau(\phi)|^2\dv<\infty.
\end{align}
Then we have
\begin{align}
\label{odd-identity-a}
\bar\nabla\bar\Delta^{2s-1}\tau(\phi)=0.
\end{align}
\end{Prop}
\begin{proof}
As in the even case the following inequality holds for an arbitrary map \(\phi\colon M\to N\) which follows by integration by parts and Young's inequality
\begin{align*}
\frac{1}{2}\int_M\eta^2|\bar\nabla^{2s+1}\bar\Delta^{s-1}\tau(\phi)|^2\dv\leq&C\int_M\eta^2|\bar\nabla^{2s-1}\bar\Delta^{s}\tau(\phi)|^2\dv 
+\frac{C}{R^2}\int_M|\bar\nabla^{2s}\bar\Delta^{s-1}\tau(\phi)|^2\dv\\
&+\int_M\eta^2|\bar\nabla^{2s}\bar\Delta^{s-1}\tau(\phi)||[\bar\Delta,\bar\nabla^{2s}]\bar\Delta^{s-1}\tau(\phi)|\dv.
\end{align*}
By iteration we obtain 
\begin{align}
\label{odd-estimate-e}
\frac{1}{2}\int_M\eta^2|&\bar\nabla^{2s+1}\bar\Delta^{s-1}\tau(\phi)|^2\dv \\
\nonumber\leq & C\int_M\eta^2|\bar\nabla\bar\Delta^{2s-1}\tau(\phi)|^2\dv  
+\frac{C}{R^2}\sum_{r=0}^{s-1}\int_M|\bar\nabla^{2s-2r}\bar\Delta^{s-1+r}\tau(\phi)|^2\dv\\
&\nonumber+\sum_{r=0}^{s-1}\int_M\eta^2|\bar\nabla^{2s-2r}\bar\Delta^{s-1+r}\tau(\phi)||[\bar\Delta,\bar\nabla^{2s-2r}]\bar\Delta^{s-1+r}\tau(\phi)|\dv,
\end{align}
which corresponds to \eqref{estimate-e} in the even case.

Now, we estimate the commutator term on the right-hand side of \eqref{odd-estimate-e} and find
\begin{align*}
\int_M\eta^2&|\bar\nabla^{2s-2r}\bar\Delta^{s-1+r}\tau(\phi)||[\bar\Delta,\bar\nabla^{2s-2r}]\bar\Delta^{s-1+r}\tau(\phi)|\dv \\
&\leq C\sum_{\sum_{l_i}+\sum_{m_j}=2s-2r}\int_M(\eta^2|\bar\nabla^{2s}\bar\Delta^{s-1}\tau(\phi)|
|\underbrace{\bar\nabla^{m_1}d\phi\star\ldots\bar\nabla^{m_{l_1}}d\phi}_{l_1-\text{times}}|\\
&\hspace{4cm}\times|\bar\nabla^{l_2}d\phi||\bar\nabla^{l_3}d\phi||\bar\nabla^{l_4+2r}\bar\Delta^{s-1}\tau(\phi)|)\dv.
\end{align*}

As a next step we calculate 
\begin{align*}
&\int_M\eta^2|\bar\nabla^{2s}\bar\Delta^{s-1}\tau(\phi)|
|\underbrace{\bar\nabla^{m_1}d\phi\star\ldots\bar\nabla^{m_{l_1}}d\phi}_{l_1-\text{times}}|
|\bar\nabla^{l_2}d\phi||\bar\nabla^{l_3}d\phi||\bar\nabla^{l_4+2r}\bar\Delta^{s-1}\tau(\phi)|\dv \\
\leq&\big(\int_M(\eta|\bar\nabla^{2s}\bar\Delta^{s-1}\tau(\phi)|)^\frac{2n}{n-2}\dv\big)^\frac{n-2}{2n}
\big(\int_M|\bar\nabla^{l_2}d\phi|^\frac{n}{l_2+1}\dv\big)^\frac{l_2+1}{n}
\big(\int_M|\bar\nabla^{l_3}d\phi|^\frac{n}{l_3+1}\dv\big)^\frac{l_3+1}{n}\\
&\times\big(\int_M(\eta|\bar\nabla^{l_4+2r}\bar\Delta^{s-1}\tau(\phi)|)^\frac{2n}{n-2(2s+1-l_4-2r)}\dv\big)^\frac{n-2(2s+1-l_4-2r)}{2n}\\
&\times\big(\int_M|\underbrace{\bar\nabla^{m_1}d\phi\star\ldots\bar\nabla^{m_{l_1}}d\phi}_{l_1-\text{times}}|^\frac{n}{\sum_{m_j}+l_1}\dv\big)^\frac{\sum_{m_j}+l_1}{n} \\
\leq &
C\big(\frac{1}{R^2}\int_M|\bar\nabla^{2s}\bar\Delta^{s-1}\tau(\phi)|^2\dv
+\int_M\eta^2|\bar\nabla^{2s+1}\bar\Delta^{s-1}\tau(\phi)|^2\dv\big)^\frac{1}{2}\\ 
&\hspace{0.4cm}\times
\big(\int_M|\bar\nabla^{l_2}d\phi|^\frac{n}{l_2+1}\dv\big)^\frac{l_2+1}{n}
\big(\int_M|\bar\nabla^{l_3}d\phi|^\frac{n}{l_3+1}\dv\big)^\frac{l_3+1}{n}\\
&\hspace{0.4cm}\times\big(
\int_M\eta^2|\bar\nabla^{2s+1}\bar\Delta^{s-1}\tau(\phi)|^2\dv 
+\sum_{q=1}^{2s+1-l_4-2r}\frac{C}{R^{2q}}\int_M|\bar\nabla^{2s+1-q}\bar\Delta^{s-1}\tau(\phi)|^2\dv\big)^\frac{1}{2}\\
&\hspace{0.4cm}\times\prod_{j=1}^{l_1}\big(\int_M|\nabla^{m_j}d\phi|^\frac{n}{m_j+l_1}\dv\big)^\frac{m_j+l_1}{n},
\end{align*}
where we first used Hölder's inequality (with \(\sum_{l_i}+\sum_{m_j}=2s-2r\)) and applied \eqref{sobolev-inequality-k} in the second step.
Note that in order to apply \eqref{sobolev-inequality-k} we have to make the assumption that \(n>4s+2\).

Combining this estimate with \eqref{odd-estimate-e} we find 
\begin{align*}
\frac{1}{2}\int_M\eta^2|&\bar\nabla^{2s+1}\bar\Delta^{s-1}\tau(\phi)|^2\dv \\
\leq & C\int_M\eta^2|\bar\nabla\bar\Delta^{2s-1}\tau(\phi)|^2\dv 
+\frac{C}{R^2}\sum_{r=0}^{s-1}\int_M|\bar\nabla^{2s-2r}\bar\Delta^{s-1+r}\tau(\phi)|^2\dv\\ 
&+C\big(\frac{1}{R^2}\int_M|\bar\nabla^{2s}\bar\Delta^{s-1}\tau(\phi)|^2\dv
+\int_M\eta^2|\bar\nabla^{2s+1}\bar\Delta^{s-1}\tau(\phi)|^2\dv\big)^\frac{1}{2}\\ 
&\hspace{0.4cm}\times
\sum_{r=0}^{s-1}\sum_{\sum_{l_i}+\sum_{m_j}=2s-2r}
\big(\int_M|\bar\nabla^{l_2}d\phi|^\frac{n}{l_2+1}\dv\big)^\frac{l_2+1}{n}
\big(\int_M|\bar\nabla^{l_3}d\phi|^\frac{n}{l_3+1}\dv\big)^\frac{l_3+1}{n}\\
&\hspace{0.4cm}\times\big(
\int_M\eta^2|\bar\nabla^{2s+1}\bar\Delta^{s-1}\tau(\phi)|^2\dv 
+\sum_{q=1}^{2s+1-l_4-2r}\frac{C}{R^{2q}}\int_M|\bar\nabla^{2s+1-q}\bar\Delta^{s-1}\tau(\phi)|^2\dv\big)^\frac{1}{2}\\
&\hspace{0.4cm}\times\prod_{j=1}^{l_1}\big(\int_M|\nabla^{m_j}d\phi|^\frac{n}{m_j+l_1}\dv\big)^\frac{m_j+l_1}{n} \\ 
\leq & C\int_M\eta^2|\bar\nabla\bar\Delta^{2s-1}\tau(\phi)|^2\dv 
+\frac{C}{R^2}\int_M|\bar\nabla^{2s}\bar\Delta^{s-1}\tau(\phi)|^2\dv\\
&+C\big(\frac{1}{R^2}\int_M|\bar\nabla^{2s}\bar\Delta^{s-1}\tau(\phi)|^2\dv
+\int_M\eta^2|\bar\nabla^{2s+1}\bar\Delta^{s-1}\tau(\phi)|^2\dv\big)^\frac{1}{2}\\
&\hspace{0.4cm}\times\epsilon
\sum_{r=0}^{s}\sum_{\sum_{l_i}+\sum_{m_j}=2s-2r} \\
&\hspace{0.4cm}\times\big(\int_M\eta^2|\bar\nabla^{2s+1}\bar\Delta^{s-1}\tau(\phi)|^2\dv 
+\sum_{q=1}^{2s+1-l_4-2r}\frac{C}{R^{2q}}\int_M|\bar\nabla^{2s+1-q}\bar\Delta^{s-1}\tau(\phi)|^2\dv\big)^\frac{1}{2},
\end{align*}
where we applied the smallness condition \eqref{odd-assumption-identity-a-1} in the second step.
Taking the limit \(R\to\infty\) while making use of the finiteness assumption \eqref{odd-assumption-identity-a-2}
and using the smallness of \(\epsilon\) we can deduce
\begin{align*}
\int_M\eta^2|\bar\nabla^{2s+1}\bar\Delta^{s-1}\tau(\phi)|^2\dv\leq C\int_M\eta^2|\bar\nabla\bar\Delta^{2s-1}\tau(\phi)|^2\dv.
\end{align*}
Combining this estimate with \eqref{odd-estimate-d} and making use of the smallness of \(\epsilon\) once
more we obtain the claim.
\end{proof}

\begin{proof}[Proof of Theorem \ref{main-theorem} (odd case)]
From this point on the proof is essentially the same as in the even case, the only difference being that one has to start the iteration
with \(j=2s-1\) using \eqref{odd-identity-a} and make use of the slightly different assumptions \eqref{odd-assumption-main-a} and \eqref{odd-assumption-main-b}.
\end{proof}

\par\medskip
\textbf{Acknowledgements:}
The author gratefully acknowledges the support of the Austrian Science Fund (FWF)
through the project P30749-N35 ``Geometric variational problems from string theory''.

\bibliographystyle{plain}
\bibliography{mybib}
\end{document}